\documentclass{article}
\usepackage{amssymb,latexsym,amsmath,amsthm,mathrsfs}
\usepackage{graphicx}
\newcommand{\comment}[1]{}

\begin{document}
\title{An easy method for finding many very large prime
numbers\footnote{Presented to the St. Petersburg Academy on March 16, 1778.
Originally published as
{\em Facillima methodus plurimos numeros primos praemagnos inveniendi},
Nova acta academiae scientiarum Petropolitanae \textbf{14} (1805),
3--10.
E718 in the Enestr{\"o}m index.
Translated from the Latin by Jordan Bell,
Department of Mathematics, University of Toronto, Toronto, Ontario, Canada.
Email: jordan.bell@gmail.com}}
\author{Leonhard Euler}
\date{}
\maketitle

1. It suffices to explain this method with a single example of numbers
contained in the formula $232aa+1$, where I shall investigate all the values
of $a$ in which this formula produces composite numbers;
with these excluded, all the remaining numbers assumed for $a$ will supply
prime numbers.

2. Now, whenever the formula $232aa+1$ yields a composite number,
an equation like this will hold:
\[
232aa+1=232xx+yy,
\]
because $232=8\cdot 29$ is a numerus idoneus.\footnote{Translator:
A number $n$ is a {\em numerus idoneus} (suitable number, convenient number,
idoneal number) if it satisfies the following criterion. Let $m$ be an odd
number relatively prime to $n$ such that for some relatively prime
$x,y$, $m=x^2+ny^2$. If the equation $m=x^2+ny^2$ has only
one solution with $x,y \geq 0$, then $m$ is a prime number.
Thus if $n$ is a numerus idoneus we can use it to find prime numbers
$m$. Gauss: a positive integer $n$ is a numerus idoneus if and only if
for every quadratic form of discriminant $-4n$, every genus consists of a single
class, see Proposition 3.24 of David A. Cox, {\em Primes of the form
$x^2+ny^2$}. It is also equivalent to $h(-4n)=2^{\mu-1}$, where $h$ is the class
number and $\mu$ is defined in the following way. Let $D \equiv 0,1 \pmod{4}$
be negative and let $r$ be the number of odd primes dividing $D$.
If $D \equiv 1 \pmod{4}$ then $\mu=r$, and if $D=-4n$, then $\mu$ is
determined as follows: $n \equiv 3 \pmod{4}$ then $\mu=r$; $n \equiv 1,2 \pmod{4}$ then $\mu=r+1$; $n \equiv 4 \pmod{8}$ then $\mu=r+1$; $n \equiv 0 \pmod{8}$
then $\mu=r+2$. See Theorem 3.22 of Cox. For more about numeri idonei see
Weil, {\em Number theory: an approach through history}, pp. 219--226.}
Therefore it will be
\[
232(aa-xx)=yy-1.
\]
Now let us put $y=1+58z$,\footnote{Translator: I don't see why Euler
says $y=1+58z$ instead of $y=1 \pm 58z$. $y=1 \pm 58z$ follows because
$58=2\cdot 29$, and both 2 and 29 are prime factors of the left hand side.}
and dividing by 232 produces
this equation: $aa-xx=\frac{1}{2}z(29z+1)$, or, because the unity needs to be taken
both positive and negative, it will be
\[
aa-xx=\frac{1}{2}z(29z \pm 1).
\]

3. Now let us take in order all values $1,2,3$ etc. for
$z$ and, because the formula $\frac{1}{2}z(29z \pm 1)$ certainly has
two factors and perhaps more,
it can be represented as a product $rs$ 
in one and perhaps more ways. Since it needs to be equal to the formula
$aa-xx=(a+x)(a-x)$, let us put $a+x=r$ and $a-x=s$, whence $a=\frac{r+s}{2}$,
and the values of $a$ to be excluded will be thus represented.
Here it is evident that both factors $r$ and $s$ need to be both even
or both odd.

4. With this noted, let us take $z=1$, and it will happen in two ways, because
of the ambiguity of the sign, as either $rs=14$ or $rs=15$. Here the first value,
since it is oddly even, is irrelevant to our purpose. On the other hand
the latter value $rs=15$ yields two exclusions, for it will either be $r=15$
and $s=1$ or $r=5$ and $s=3$, from which the values 8 and 4 need to be excluded
for $a$.

5. Now let $z=2$ and it will be $rs=57$ or $rs=59$. The first admits two
resolutions:

\begin{center}
\begin{tabular}{cc}
$r=57$,&$r=19$,\\
$s=1$,&$s=3$,
\end{tabular}
\end{center}

the latter one: $r=59$ and $s=1$, from which these three exclusions arise:
$a=29, a=11, a=30$.

6. Now let $z=3$, and it will be either $rs=3\cdot 43$ or $rs=3\cdot 44$,
whose resolutions happen like this:

\begin{center}
\begin{tabular}{cc|cc}
$r=129$,&$r=43$,&$r=66$,&$r=22$,\\
$s=1$,&$s=3$,&$s=2$,&$s=6$,
\end{tabular}
\end{center}

from which these four exclusions follow

\[
a=65, \quad a=23, \quad a=34, \quad a=14.
\] 

7. Let $z=4$, and it will be either $rs=230$ or $rs=234$. This numbers,
since they are oddly even, give no exclusions.

8. Now let $z=5$, and it will be $rs=5\cdot 72$ or $rs=5\cdot 73$.
The first of these yields the following resolutions:

\begin{center}
\begin{tabular}{cccccc}
$r=180$,&$r=90$,&$r=60$,&$r=36$,&$r=30$,&$r=20$,\\
$s=2$,&$s=4$,&$s=6$,&$s=10$,&$s=12$,&$s=18$.
\end{tabular}
\end{center}

Then the values of $a$ that need to be excluded are therefore $91,47,33,23,21,19$.
The other value gives $r=73,s=5$ or $r=365,s=1$,
from which the exclusions 39 and 183 arise.

8a.\footnote{Translator: The original has two sections 8.} Let $z=6$,
and it will be either $rs=3\cdot 173$ or $rs=3\cdot 175$,
whence the following resolutions:

\begin{center}
\begin{tabular}{cc|cccccc}
$r=519$,&173,&$r=525$,&175,&105,&75,&35,&25,\\
$s=1$,&3,&$s=1$,&3,&5,&7,&15,&21.
\end{tabular}
\end{center}

Therefore these exclusions arise:

\[
260,\qquad 88, \qquad 263, \qquad 89, \qquad 55, \qquad 41, \qquad 25, \qquad
23.
\]

9. Now let $z=7$, and it will be either $rs=7\cdot 101$ or $rs=7\cdot 102$.
Here the first case alone gives these resolutions: 

\begin{center}
\begin{tabular}{cc}
$r=707$,&$r=101$,\\
$s=1$,&$s=7$,
\end{tabular}
\end{center}

and then the exclusions will be 354 and 54. So that we can prescribe some bounds
for ourselves, we will omit in the sequel all exclusions greater than 300.

10. Now let $z=8$, and it will be either $rs=4\cdot 3\cdot 7\cdot 11$ or
$rs=4\cdot 233$.
Then the following resolutions arise:

\begin{center}
\begin{tabular}{cccc|c}
$r=462$,&154,&66,&42,&$r=466$,\\
$s=2$,&6,&14,&22,&$s=2,$
\end{tabular}
\end{center}

from which we conclude the following exclusions:

\[
a=232, \qquad 80, \qquad 40, \qquad 32, \qquad 234.
\]

11. Let $z=9$, and it will be either $rs=9\cdot 130$ or $rs=9\cdot 131$.
Then for the second

\begin{center}
\begin{tabular}{cc}
$r=393$,&131,\\
$s=3$,&9,
\end{tabular}
\end{center}

one therefore excludes $a=198,70$.

12. Let $z=10$, and it will be either $rs=3\cdot 5\cdot 97$ or
$rs=5\cdot 17\cdot 17$, then

\begin{center}
\begin{tabular}{ccc|cc}
$r=485$,&291,&97,&$r=289$,&85\\
$s=3$,&5,&15,&$s=5$,&17, 
\end{tabular}
\end{center}

whence one excludes
\[
a=244,\qquad 148,\qquad 56,\qquad 147,\qquad 51.
\]

13. Let $z=11$, and it will be either $rs=3\cdot 11\cdot 53$ or $rs=32\cdot 5\cdot 11$, then

\begin{center}
\begin{tabular}{ccc|ccccccc}
$r=583$,&159,&53,&440,&220,&176,&88,&80,&110,&44,\\
$s=3$,&11,&33,&4,&8,&10,&20,&22,&16,&40,
\end{tabular}
\end{center}

hence

\begin{center}
\begin{tabular}{cccccccccc}
$a=293$,&85,&43,&222,&114,&93,&54,&51,&63,&42.
\end{tabular}
\end{center}

14. Let $z=12$, and it will be either $rs=6\cdot 347$ or
$rs=6\cdot 349$, neither of which gives any exclusion.

15. Let $z=13$, and it will be either $rs=13\cdot 4\cdot 47$ or $rs=27\cdot 13\cdot 7$, then

\begin{center}
\begin{tabular}{c|cccccc}
$r=94$,&351,&273,&189,&117,&91,&63,\\
$s=26$,&7,&9,&13,&21,&27,&39;
\end{tabular}
\end{center}

one therefore excludes

\begin{center}
\begin{tabular}{ccccccc}
60,&179,&141,&101,&69,&59,&51.
\end{tabular}
\end{center}

16. The calculation can be continued on in this way as far as it seems
fit. We shall represent the exclusions arising from all the values of $z$
in the following table:

\begin{center}
\begin{tabular}{r|p{10cm}}
$z$&Exclusions\\
\hline
1&4, 8,\\
2&11, 29, 30,\\
3&14, 23, 34, 65,\\
4&$-\,  -\,  -\,  -$\\
5&19, 21, 23, 33, 39, 47, 91, 183,\\
6&23, 25, 41, 55, 88, 89, 260, 263,\\
7&54,\\
8&32, 40, 80, 232, 234,\\
9&70, 198,\\
10&51, 56, 147, 148, 244,\\
11&42, 43, 51, 54, 63, 85, 93, 114, 222, 293,\\
13&51, 59, 60, 69, 101, 141, 179,\\
14& 54, [57], 58, 66, 78, 102, 135, 162, 206, 207, [286],\\
15&64, 68, 88, 116, 236,\\
16&61, 77, 103, 161, 191,\\
17&120, 132, 168,\\
18&265, 266,\\
19&75, 80, 88, [117], 147, 243,\\
21&80, 83, 85, 97, 103, 109, 128, 145, 163, 221, 235, 272,\\
22&84, [85], 96, 123, 276,\\
23&90, 122, 178,\\
24&92, 140, 154,\\
25&98, 110, 154, 194,\\
26&108, 145,\\
27&103, 105, 111, 125, 129, 147, [165], 203, 209, 241,\\
29&122, 134, 166, 218, 225,\\
30&122, 131, 133, 146, 187, 214,\\
31&240,\\
32&139, 224,\\
33&256,\\
34&130, 141, 190, [202],\\
35&137, 143, 162, 169, 247, 271,\\
37&141, 145, 171, 287,
\end{tabular}
\end{center}

\begin{center}
\begin{tabular}{r|p{10cm}}
$z$&Exclusions\\
\hline
38&212,\\
39&154, 202,\\
40&156, 174, 178, 242,\\
41&162, 186,\\
42&160, 220, 268,\\
43&164, [181], 195, 199, 211, 284,\\
45&208,\\
46&191, 217, 257,\\
47&184,\\
48&185, 241, 253,\\
49&[192], 260,\\
50&191, 193, [215, 225, 279, 297],\\
51&196, 236, 292,\\
53&202, 223, 245, 260,\\
54&224,\\
56&216, 240,\\
58&224, 236, 296,\\
59&225, 231, 233, 273,\\
61&239, 241, 282,\\
62&237,\\
63&240,\\
64&246, 282,\\
65&252,\\
66&274,\\
67&282,\\
69&263, 265, 295, 263, 265,\\
71&278,\\
73&286,\\
75&286, 290,\\
78&298.
\end{tabular}
\end{center}

17. Among the numbers to be excluded for the letter $a$, some occur
two or three times. This happens whenever the number $232aa+1$
has three or more factors; thus we see that the number 265 is excluded
first by the value $z=18$, then twice by the value $z=69$.
In fact in this case it turns out that the number
$232\cdot 265^2+1$
is made from the three factors
$59\cdot 461\cdot 599$,
and even aside from this form the following three 
are recovered:

\begin{center}
\begin{tabular}{ll}
$1\,^{\circ}$&$232\cdot 256^2+1043^2$,\\
$2\,^{\circ}$&$232\cdot 34^2+4003^2$,\\
$3\,^{\circ}$&$232\cdot 35^2+4001^2,$
\end{tabular}
\end{center}

18. Let us also record in a table all the numbers that are to be excluded
according to the following arrangement, and we note with an asterisk all
those with occur more than once:

\begin{center}
\begin{tabular}{r|l}
0&4, 8,\\
1&11, 14, 19,\\
2&21, 23*, 25, [29],\\
3&30, 32, 33, 34, 39,\\
4&40, 41, 42, 43, 47,\\
5&51*, 54*, 55, 56, [57], 58, 59,\\
6&60, 61, 63, 64, 65, 66, 68, 69,\\
7&70, 75, 77, 78,\\
8&80*, 83, 84, 85*, 88*, 89,\\
9&90, 91, 92, 93, 96, 97, 98,\\
10&101, 102, 103*, 105, 108, 109,\\
11&110, 111, 114, 116, [117],\\
12&120, 122*, 123, 125, 128, 129,\\
13&130, 131, 132, 133, 134, 135, 137, 139,\\
14&140, 141*, 143, 145*, 146, 147*, 148,\\
15&154*, 156,\\
16&160, 161, 162, 163, 164, [165], 166, 168, 169,\\
17&171, 174, 178*, 179,\\
18&[181], 183, 184, 185, 186, 187,\\
19&190, 191*, [192], 193, 194, 195, 196, 198, 199,\\
20&202*, 203, 206, 207, 208, 209,\\
21&211, 212, 214, [215], 216, 217, 218,\\
22&220, 221, 222, 223, 224*, 225*,\\
23&231, 232, 233, 234, 235, 236*, 237, 239,\\
24&240*, 241*, 242, 243, 244, 245, 246, 247,\\
25&252, 253, 256, 257,\\
26&260*, 263*, 265*, 266, 268,\\
27&271, 272, [273], 274, 276, 278, [279],\\
28&282*, 284, 286*, 287,\\
29&290, 292, 293, 295, 296, [297], 298.
\end{tabular}
\end{center}

Therefore with these numbers excluded, all the remaining written in place of $a$
in the formula $232aa+1$ produce prime numbers. The following table
thus exhibits these values of $a$:

\begin{center}
\begin{tabular}{rrrrrrrrrrr}
1,&2,&3,&5,&6,&7,&9,&10,&12,&13,&15,\\
16,&17,&18,&20,&22,&24,&26,&27,&28,&31,&35,\\
36,&37,&38,&44,&45,&46,&48,&49,&50,&52,&53,\\
62,&67,&71,&72,&73,&74,&76,&79,&81,&82,&86,\\
87,&94,&95,&99,&100,&104,&106,&107,&112,&113,&115,\\
118,&119,&121,&124,&126,&127,&136,&138,&142,&144,&149,\\
150,&151,&152,&153,&155,&157,&158,&159,&167,&170,&172,\\
173,&175,&176,&177,&180,&182,&188,&189,&197,&200,&201,\\
204,&205,&210,&213,&219,&226,&227,&228,&229,&230,&238,\\
248,&249,&250,&251,&254,&255,&258,&259,&261,&262,&264,\\
267,&269,&270,&275,&277,&280,&281,&283,&285,&288,&289,\\
291,&294,&299;&&&&&&&&
\end{tabular}
\end{center}

the last of these numbers comes to about 20 million.\footnote{Translator: If
$a=299$ then $232aa+1=$ 20,741,033.}

\end{document}